\documentclass[12pt]{article}
\usepackage{latexsym}
\usepackage{amssymb, amsmath, xspace, lscape, latexsym, mathrsfs, epsf}
\usepackage{overpic}
\usepackage{hyperref}
\newcommand{\De}{\Delta}

\newcommand{\pr}{\prime}
\newcommand{\non}{\nonumber}
\newcommand{\bea}{\begin{eqnarray}}
\newcommand{\eea}{\end{eqnarray}}
\def\beq#1#2\eeq{
        \begin{equation}
        \label{#1}
            #2
        \end{equation}}

\newcommand{\al}{\alpha}
\newcommand{\bt}{\beta}
\newcommand{\te}{\theta}

\newcommand{\Ga}{\Gamma}

\newcommand{\rex}{{\rm e}}

\newcommand{\Int}{\int_{0}^{1}}

\newcommand{\bq}{\begin{eqnarray}}
\newcommand{\eq}{\end{eqnarray}}
\newcommand{\nn}{\nonumber}
\newcommand{\ba}{\begin{array}}
\newcommand{\ea}{\end{array}}

\newcommand{\tJ}{\textsf{J}}

\def\btheor#1\etheor{
        \begin{theor}
            #1
        \end{theor}
    }

    \def\bsled#1\esled{
        \begin{sled}
            #1
        \end{sled}   }

\newtheorem{theorem}{Theorem}

\newtheorem{lemma}{Lemma}
\newtheorem{prop}{Proposition}

\def\hm#1{#1\nobreak\discretionary{}{\hbox{\m@th$#1$}}{}}
\def\mi#1{\discretionary{\hbox{\m@th$#1$}}{\hbox{\m@th$#1$}}{}}

\vspace{4ex}

\begin{document}

\title{\bf Painlev\'e VI and the Unitary Jacobi ensembles}
\author{
Yang Chen\\
\small{\textit{Department of Mathematics, Imperial College London,}} \\
\small{\textit{180 Queen's Gates, London SW7 2BZ, UK}}\\
\small{ychen@imperial.ac.uk}\\
\vspace{1mm}
and
\\
\vspace{1mm}
Lun Zhang \\
\small{\textit{Department of Mathematics, Katholieke Universiteit Leuven,}} \\
\small{\textit{Celestijnenlaan 200 B, 3001 Leuven, Belgium}}\\
\small{lun.zhang@wis.kuleuven.be} } \maketitle

\begin{abstract}
The six Painlev\'e transcendants which originally appeared in the
studies of ordinary differential equations have been found numerous
applications in physical problems. The well-known examples among
which include symmetry reduction of the Ernst equation which arises
from stationary axial symmetric Einstein manifold and the spin-spin
correlation functions of the two-dimensional Ising model in the work
of McCoy, Tracy and Wu.

The problem we study in this paper originates from random matrix theory,
namely, the smallest eigenvalues distribution of the finite $n$
 Jacobi unitary ensembles which was first investigated by Tracy and
Widom. This is equivalent to the computation of the probability that
the spectrum is free of eigenvalues on the interval $[0,t]$. Such
ensembles also appears in multivariate statistics known as the
double-Wishart distribution.

We consider a more general model where the Jacobi weight is
perturbed by a discontinuous factor and study the associated finite
Hankel determinant. It is shown that the logarithmic derivative of
Hankel determinant satisfies a particular $\sigma-$form of
Painlev\'e VI, which holds for the gap probability as well. We also
compute exactly the leading term of the gap probability as $t\to
1^-$.

\end{abstract}

\vfill\eject
\section{Introduction and statement of results}
The discovery by Tracy and Widom \cite{TWA}, \cite{TWB} in the
1990's of probability laws that govern the extreme eigenvalues of
certain Hermitian random matrices triggered a renewed interest in
random matrix theory and has lead to an explosion of activity in
recent years; cf. \cite{chen+bai} for a recent account regarding the
application of random matrix theory to wireless communications and
multi-variate statistics. In particular, Painlev\'e transcendants
play an important role in these studies.

Painlev\'e transcendants refer to solutions of six nonlinear second
order ordinary differential equations of the form
$$
y''=R(y,y',t).
$$
Here $'$ denotes $\frac{d}{dt}$ and $R$ is a function rational in
$y,$ $y'$ and analytic in $t.$ All Painlev\'e transcendants enjoy
the celebrated Painlev\'e property which asserts the only movable
singularities are poles.

There are broadly speaking two well-known physical applications of
the Painlev\'e equations. One is the reduction of the integrable
Ernst-Weyl equations to Painlev\'e VI by Cosgrove \cite{cos} and the
reduction of the integrable Loewner-Konopelchenko-Rogers system
\cite{LKR} to a matrix form of Painlev\'e VI by Schief
\cite{Schief}. For a review on the reduction of the Ernst equation,
we refer to \cite{Mason}. The other application is to the
two-dimensional Ising model where McCoy, Tracy and Wu \cite{MTW}
showed that the spin-spin correlation function expressed as the
Fredholm determinant of certain integral operator is a particular
Painlev\'e III.

In this paper, we are also concerned with Fredholm determinant of a certain
finite rank integral operator, which is equivalent to the determinant of
a Hankel matrix generated by a weight of the form
\begin{equation}\label{weight function}
x^{\al}(1-x)^{\bt}(A+B\theta(x-t)),\quad x\in[0,1]
\end{equation}
where $t\in[0,1]$, $\alpha>0$, $\beta>0$ and $\te(\cdot)$ denotes
the Heaviside step function.

We take the shifted Jacobi weight,
$x^{\al}(1-x)^{\bt},$ to be  the ``reference'' weight $w_0$ and call
$w_0w_\tJ$ a perturbed Jacobi weight, where
$w_{\tJ}(x;t)=A+B\te(x-t)$ is the ``jump'' function. The constants
$A$ and $B$ are chosen in such a manner that $A\geq 0$ and $A+B \geq
0$, which ensures the real Hankel form
\begin{equation*}
\sum_{i,j=0}^{n-1}c_i\:c_j\:\mu_{i+j}(t), \quad c_i \in
\mathbb{R},\quad  i=0,1,...,
\end{equation*} is positive definite
and the Hankel matrix defined below
\begin{equation}
{\cal H}_n:=\left(\mu_{i+j}(t)\right)_{i,j=0}^{n-1}\:,\nn
\end{equation}
is invertible, where
\begin{equation}\label{moments}
\mu_i(t):=\Int\:x^{\al+i}(1-x)^{\bt}[A+B\te(x-t)]dx,
\end{equation}
for $i=0,1,2,...$. An evaluation of the integrals in \eqref{moments}
gives
\begin{equation}
\mu_i(t)=A\Ga(i+\al+1)+\frac{B}{1+\bt}(1-t)^{1+\bt}\:t^{i+\al}\:
_2F_1(-\al-i,1;2+\bt;1-1/t)\nn,
\end{equation}
where $\Gamma(x)$ denotes the Gamma function and $_2F_1$ is the
hypergeometric function; see \cite{Abr+Ste}.

It is a well-known fact that the Hankel determinant
\begin{equation}\label{Hankel det}
D_n(t)=\det{\cal H}_n=\det\left(\mu_{i+j}(t)\right)_{i,j=0}^{n-1}
\end{equation}
has the following multiple integral representation \cite{Szego}
\begin{equation}\label{integral representation}
\frac{1}{n!}\int_{[0,1]^n}
[\Delta_n(x)]^2\prod_{k=1}^{n}x_k^{\al}(1-x_k)^{\bt}w_{\tJ}(x_k,t)dx_k,
\end{equation}
where \bea \Delta_n(x)=\prod_{1\leq i<j\leq n}(x_j-x_i).\nn
\eea Hence, we can also view $D_n(t)$ as the partition function for
the random matrix ensemble with the joint eigenvalue distribution
$$
\prod_{1\leq j<k\leq
n}(x_j-x_k)^2\prod_{l=1}^{n}x_l^{\al}(1-x_l)^{\bt}\:(A+B\theta(x_l-t))dx_l,
$$
which is of interest in statistical physics.

Let
\begin{equation}\label{Hn}
H_n(t):=t(t-1)\frac{d}{dt}\ln D_n(t).
\end{equation}
It is the aim of this paper to produce a second order differential
equation (non-linear) that is satisfied by the above logarithmic
derivative of $D_n(t)$. Our main result is the following:
\begin{theorem}\label{thm of sigmaform}
The $\sigma$ function, defined as a translation in $t$ of $H_n,$
\begin{align} \label{htil-def}
\sigma(t):=H_n(t) +d_1 t+d_2
\end{align}
with
\begin{align*}
d_1 &=  - n (n + \al + \bt )-\frac{(\al + \bt)^2}{4},  \\
d_2 & = \frac{1}{4} \biggl[ 2n (n + \al + \bt ) + \bt (\al + \bt)
\biggr],
\end{align*}
satisfies the following Jimbo-Miwa-Okamoto (\cite{Jimbo-Miwa},
\cite{Okamoto}) $\sigma-$form of Painlev\'e VI
\begin{align}
& \sigma'\biggl\{t (t-1) \sigma''\biggr\}^2 + \biggl\{ 2 \sigma'
\left( t \sigma' - \sigma \right) - \sigma'^2 - \nu_1 \nu_2 \nu_3
\nu_4 \biggr\}^2
\nonumber\\
& = \left(\sigma' + \nu_1^2 \right)\left(\sigma' + \nu_2^2
\right)\left(\sigma' + \nu_3^2 \right)\left(\sigma' + \nu_4^2
\right) \label{sigmaform}
\end{align}
with
\begin{equation*}
\nu_1 = \frac{\al + \bt}{2}, \quad \nu_2 = \frac{\bt - \al}{2},
\quad \nu_3 = \nu_4= \frac{2n + \al + \bt }{2}
\end{equation*}
and the initial conditions \bea
\sigma(0)=d_2,\quad\sigma'(0)=d_1.\nn \eea
\end{theorem}
Note that the constants $A$ and $B$ do not appear in the
differential equation \eqref{sigmaform}.

In the situation where $A=0$ and $B=1$, the quotient $ D_n(t)/D_n(0)
$ stands for the gap probability that the interval $[0,t]$ is free
of eigenvalues or, equivalently, the least eigenvalue of (the finite
$n$) Jacobi unitary ensembles is $\geq t.$
By a fundamental result of Gaudin and Mehta \cite{Mehta}, the gap
probability can be formulated as the following Fredholm determinant:
\begin{equation}
\frac{D_n(t)}{D_n(0)}=\det(I-K_{[0,t]}^{(n)}).\nn
\end{equation}
The integral operator $K_{[0,t]}^{(n)}$ acts on functions $f\in
L^{2}(0,t)$ as
\begin{equation} K_{[0,t]}^{(n)}f(x)=\int_{0}^{t}a_n
\frac{\phi_{n}(x)\phi_{n-1}(y)-\phi_{n-1}(x)\phi_{n}(y)}{x-y}f(y)dy,\nn
\end{equation}
where $a_n$ is a constant explicitly known and $\phi_n(x)$
is a constant multiple of
$x^{\al/2}\:(1-x)^{\bt/2}\:P_n^{(\al,\bt)}(x),$ with
$P_n^{(\al,\bt)}(x)$ being the Jacobi polynomials. Our next result
deals with the asymptotic behavior of gap probability as $t\to 1^-$:
\begin{theorem}\label{asy of gap probability}
Suppose that $A=0$ and $B=1$, we have
\begin{align}\label{asy near 1}
\frac{D_n(t)}{D_n(0)}\sim
2^{2n\al}\:\frac{F(0,\bt)K(0,\bt,n)}{F(\al,\bt)K(\al,\bt,n)}(1-t)^{n(n+\bt)},
\end{align}
as $t\to 1^-$, where
\begin{align}
F(\al,\bt)&=\frac{\Ga((\al+\bt+1)/2)G^2((\al+\bt+1)/2)G^2(1+(\al+\bt)/2)}
{G(\al+\bt+1)G(\al+1)G(\bt+1)},\label{F(al,bt)}\\
K(\al,\bt,n)&=\frac{G(n+1)G(n+\al+1)G(n+\bt+1)G(n+\al+\bt+1)}
{G^2(n+(\al+\bt+1)/2)G^2(n+1+(\al+\bt)/2) \Ga(n+(\al+\bt+1)/2)}.
\label{K(al,bt)}
\end{align}
Here, $G(x)$ denotes the Barnes G-function:
\begin{equation*}
G(x+1)=(2\pi)^{x/2}e^{-(x(x+1)+\gamma x^2)/2}\prod_{n=1}^{\infty}
\Big[\Big(1+\frac{x}{n}\Big)^n e^{-x+x^2/(2n)}\Big],
\end{equation*}
where $\gamma$ is the Euler-Mascheroni constant.
\end{theorem}

We mention here that the so-called Dyson-Widom constant; the coefficient of $(1-t)^{n(n+\beta)}$ of
(1.8) was conjectured
in \cite{Oleg} for the gap problem, but with the more general hypergeometric
kernel. However that result cannot be applied here since the 
parameters there and those of our problem are different.

The differential equation satisfied by the logarithmic derivative of
$D_n(t)$ when $A=0$ and $B=1$ is one of the problems first
investigated by Tracy and Widom in \cite{TWFred} using operator
approach of \cite{TWA} and \cite{TWB}. This method brings in the
resolvent of $K_{[0,t]}^{(n)}$ and distributional kernel, gives rise
to a third order ordinary differential equation. It was not clear
how a first integral could be found, as one expects to find a second
order equation. Recently, a first integral was found in \cite{Oleg}
for the Gap problem for hypergeometric kernel.

Later, Haine and Semengue \cite{h+s}, made use of the
``infinitely-many-times'' approach of Adler, Shiota and van Moerbeke
\cite{Adl+Van}, popular amongst the experts in integrable systems,
to investigate the gap problem in the Jacobi polynomials ensemble.
This approach involves the notion of vertex operators, an important
tool in Sato's theory. Another third order ordinary differential
equation was found. Eliminating the third derivatives from both
equations, Haine and Semengue, found that the resulting second order
ordinary differential equation is a particular Painlev\'e VI.

The appearance of Painlev\'e VI could perhaps be guessed at based on
the following considerations: In the original or ``unperturbed''
problem, where $B=0,$ the Jacobi polynomials satisfies a second
order ordinary differential equation with regular singular points at
$0,$ $1,$ and $\infty.$ It is to be expected that ``switching'' on
the discontinuity ($B>0$), an extra regular singular point is
introduced at $t$. Because the resulting ordinary differential
equation has four regular singular points, one of which moving,
according to the isomonodromy theory of Jimbo and Miwa
\cite{Jimbo-Miwa}, one expects to find Painlev\'e VI and allied
functions appearing in the coefficients of the ordinary differential
equation. However, it is not immediately clear how the correct
parameters may be identified, and it is also not clear which
quantities related to the orthogonal polynomials is that particular
Painlev\'e VI.

Because of the discontinuity of the weight at $t$ and the $t$
dependence of the moment sequence $\{\mu_i(t)\}$, $i=1,2,...,$ the
coefficients of polynomials $\{P_n(z)\}$ orthogonal with respect to
$w_0w_\textsf{J}$ also depend on $t.$ We normalize our monic
$P_n(z)$ to be
\begin{equation}
P_n(z)=z^{n}+\textsf{p}_1(n,t)z^{n-1}+...+P_n(0).\nn
\end{equation}
It turns out that
\begin{equation}
H_n(t)= t(t-1)\frac{d}{dt}\ln
D_n(t)=-(2n+\al+\bt)\textsf{p}_1(n,t)-n(n+\al);\nn
\end{equation}
see Theorem \ref{H_n expression} in Section \ref{Ladder operators}.
This, together with \eqref{htil-def}, implies that our
$\sigma-$function is up to a translation in $t$ essentially
$\textsf{p}_1(n,t)$, the coefficient of $z^{n-1}$ of monic
orthogonal polynomials $P_n(z)$. Comparing with the Tracy-Widom
approach, we see that their resolvent evaluated at the same point is
also interpreted here as $\textsf{p}_1(n,t)$.

The approach in this paper is based on the ladder operators of
orthogonal polynomials and the associated compatibility conditions
$(S_1),$ $(S_2)$ and $(S_2')$, the detailed derivation of which can
be found in \cite{basor+chen}. We find through the ladder operators
there are four auxiliary quantities $x_n(t),\;y_n(t),\;R_n(t)$ and
$r_n(t).$ This appears to create complications. But the
compatibility conditions, as we shall see later, ultimately provides
a determination of these through a system of difference equations.

For a comparison between the ladder operator theory and the
isomonodromy theory of Jimbo and Miwa \cite{Jimbo-Miwa} in
 a particular Hermitian random ensembles, see
\cite{chen+its}. There the polynomials are orthogonal with respect
to a ``singularly'' deformed Laguerre weight;
$x^{\al}\rex^{-x}\:\rex^{-s/x},$ over $\mathbb{R}_{+},$ and one
finds a Painlev\'e III. The reader is also referred to \cite{igor}
for a comparison between the Tracy-Widom theory and that of Adler,
Shiota and van Moerebeke.

In the development of what follows, we would
like to remind the reader to keep in mind that the recurrence
coefficients $\alpha_n(t)$ and $\beta_n(t)$ (see \eqref{recurrence
relation} below) play a central role and attention should be focused
on expressing these in terms of the auxiliary quantities mentioned
above. The reader will see that the auxiliary quantities appear
naturally in the ladder operator approach. We will show that one of
the auxiliary quantities $x_n(t)$, intimately related to the
recurrence coefficients, is up to a linear fractional transformation
a particular Painlev\'e VI.

This rest of the paper is organized as follows. In the next section,
we introduce the notations to be used throughout the paper and
derive two Toda-type equations; these are differential-difference
equations involving $\al_n(t)$ and $\bt_n(t)$. In Section
\ref{Ladder operators}, a system of non-linear difference equations
involving the auxiliary quantities are found and are used to obtain
expressions of the recurrence coefficients amongst other things.
Based on these facts, we prove Theorem \ref{thm of sigmaform} and
show that $x_n(t)$--one of the auxiliary quantities--is actually a
particular Painlev\'e VI in Section \ref{derivation of PVI}. We
conclude this paper in Section \ref{proof of thm 2} with the
proof of Theorem \ref{asy of gap probability}.

\setcounter{equation}{0}

\section{Preliminaries: notations and time evolution}
The purpose of this section is to derive two universal Toda-type
equations which involve the recurrence coefficients and the
auxiliary quantities $R_n(t)$ and $r_n(t)$. We keep $w_0$ quite
general.

The general Toda-hierarchy, obtained through the multiplication of
$w_0$ by $ \exp\left(-\sum_{j=0}^{\infty}t_j\:x^{j}\right) $ can be
found for example, in \cite{Moser}, \cite{Adl+Van} and \cite{h+h}.
See also \cite{h+s} for a discussion of the approach of
\cite{Adl+Van} in relation to Sato's theory.

To begin we consider polynomials $\{P_{i}(x)\}$ orthogonal with
respect to the weight
$$w_0(x)(A+B\theta(x-t))$$
 on [0,1]. The
weight $w_{0}$ will be known as the ``reference'' weight, and
$$
w_{\textsf{J}}(x;t):=A+B\theta(x-t)
$$
as the ``jump'' factor.

The orthogonality condition is
\begin{equation}
\Int\:P_i(x)P_j(x)w_0(x)w_{\textsf{J}}(x;t)dx=h_i(t)\delta_{i,j},\quad
h_i(t)>0
\end{equation}
for $i=0,1,...$, and the $t$ dependence through $w_{\textsf{J}},$
induces $t$ dependence on the coefficients. We normalize our monic
polynomials as
\bea
P_n(z)=z^n+\textsf{p}_{1}(n,t)z^{n-1}+...+P_n(0),
\eea
although we
do not always display the $t$ dependence of coefficients of
$z^{n-1}$.

An immediate consequence of the orthogonality condition is the three
terms recurrence relations;
\begin{equation}\label{recurrence relation}
zP_{n}(z)=P_{n+1}(z)+\al_nP_n(z)+\bt_{n}P_{n-1}(z),
 \end{equation} with the
initial conditions \bea P_0(z)=1,\qquad \bt_0P_{-1}(z)=0. \eea An
easy consequence of the recurrence relation is
\begin{equation}\label{al_n-p_n}
\al_n(t)=\textsf{p}_1(n,t)-\textsf{p}_1(n+1,t),
\end{equation} and a
telescopic sum of the above equation (bearing in mind that
$\textsf{p}_1(0,t)=0),$ leaves
\begin{equation}\label{sum al_j}
-\sum_{j=0}^{n-1}\al_j(t)=\textsf{p}_1(n,t).
\end{equation}
First let us discuss the derivatives of $\al_n$ and $\bt_n$ with
respect to $t,$ as this yields the simplest equations. We keep $w_0$
quite general, as long as the moments, \bea
\mu_i(t):=\int_{0}^{1}x^{i}w_0(x)w_{\textsf{J}}(x;t)dx,\quad
i=0,1,..., \eea of all orders exist. Taking a derivative of $h_n$
with respect to $t$ gives, \bea h_n'(t)=-Bw_0(t)P_n(t,t)^2 \eea or
\begin{equation} \label{log derivative of hn}
(\log h_n)'=-R_n,
\end{equation}
where
\begin{equation}
R_n(t):=B w_0(t)\frac{P_n(t,t)^2}{h_n(t)}.
\end{equation}
Because $\bt_n=h_n/h_{n-1}$, we have the first Toda equation, \bea
\bt_n'=(R_{n-1}-R_n)\bt_n.\nn \eea

Let $D_{n}(t)$ be the Hankel determinant given in \eqref{Hankel
det}. It is well-known that
$$
D_{n} = \prod_{ i= 0}^{n-1} h_{i}(t).
$$
In view of \eqref{log derivative of hn}, it is easily seen that
\begin{equation}\label{d log D_n}
\frac{d}{dt}\log D_n(t)=-\sum_{j=0}^{n-1}R_j.
\end{equation}
Also,
\begin{align}
0&=\frac{d}{dt}\Int\:P_n(x)P_{n-1}(x)w_0(x)w_{\textsf{J}}(x;t)dx\nonumber\\
&=-Bw_0(t)P_n(t,t)P_{n-1}(t,t)
+h_{n-1}\;\frac{d}{dt}\textsf{p}_1(n,t)\nonumber
\end{align} and therefore
\bea
\frac{d}{dt}\textsf{p}_1(n,t)=Bw_0(t)\frac{P_n(t,t)P_{n-1}(t,t)}{h_{n-1}(t)}=:r_n(t).\nn
\eea This, together with \eqref{al_n-p_n}, implies
\begin{equation}
\frac{d}{dt}\al_n=r_{n}-r_{n+1}.\nn
\end{equation}
To summarize we have
\begin{theorem}
The recurrence coefficients $\al_n,\;\bt_n$ and the auxiliary
quantities $R_n,$$r_n$ satisfy the coupled Toda equations,
\begin{align}\label{differential equ for btn}
\bt_n'&=(R_{n-1}-R_n)\bt_n,
\\
\label{derivative of aln } \al_n'&=r_n-r_{n+1},
\end{align}
and
\begin{equation}
\label{derivative of p_1n} \textsf{p}_1'(n,t)=r_n(t).
\end{equation}
\end{theorem}
Note that the above system \eqref{differential equ for
btn}--\eqref{derivative of aln } does not close, unlike the Toda
equations under the first time flow.

It will transpire from next section that we will encounter two more
auxiliary quantities $x_n$ and $y_n.$ The ultimate aim of the next
section are the expressions of $\al_n$ and $\bt_n$ in terms of the
auxiliary quantities, and certain identities involving $R_n$, $r_n$,
$x_n$ and $y_n$. \setcounter{equation}{0}

\section{Ladder operators, compatibility conditions and difference equations}
\label{Ladder operators}

In this section, we give an account for a recursive algorithm for
the determination of the recurrence coefficients $\al_n$ and
$\bt_n$, based on a pair of ladder operators and the associated
supplementary conditions; see $(S_1),\:(S_2)$ and $(S_2')$ below.
Such operators have been derived by various authors over many years.
For a quick guide to the relevant literature, the reader is directed
to the references of \cite{chen+its} and \cite{bacheh}. In fact
Magnus \cite{Mag1} traced this back to Laguerre. The derivation
where the weight has discontinuities can be found in
\cite{basor+chen}.

We find the form of the ladder operators and the compatibility conditions set out below
well suited for our purpose.

Let $w(x)=w_0(x)w_{\textsf{J}}(x;t)$. For a sufficiently
well-behaved ``reference'' weight $w_0(x)$ (see \cite{chen2} for a
precise statement) of the form
$$
w_0(x)={\rm e}^{-\textsf{v}_0(x)},
$$
the lowering and raising
operators are
\begin{align}
P_n'(z)&=-B_n(z)P_n(z)+\bt_nA_n(z)P_{n-1}(z),\\
P_{n-1}'(z)&=[B_n(z)+\textsf{v}_0^{\pr}(z)]P_{n-1}(z)-A_{n-1}(z)P_n(z),
\end{align} with
\begin{align}
A_n(z)&:=\frac{R_n(t)}{z-t}+\frac{1}{h_n}\Int\:
\frac{\textsf{v}_0'(z)-\textsf{v}_0'(y)}{z-y}
P_n^2(y)w(y)dy,\\
B_n(z)&:=
\frac{r_n(t)}{z-t}+\frac{1}{h_{n-1}}\Int\:\frac{\textsf{v}_0'(z)-\textsf{v}_0'(y)}{z-y}
P_n(y)P_{n-1}(y)w(y)dy,
\end{align}
where we have assumed that $w_0(0)=w_0(1)=0.$ Additional terms would
have to be included in the definitions of $A_n(z)$ and $B_n(z)$, if
$w_0(0), w_0(1)\neq 0.$
Due to the discontinuity of $w_{\textsf{J}}$ at $t$, $A_n(z)$ and
$B_n(z)$ have apparent poles at $t$ with residues $R_n(t)$ and
$r_n(t)$, respectively.

A direct calculation produces two fundamental supplementary
(compatibility) conditions valid for all $z;$
$$
B_{n+1}(z)+B_n(z)=(z-\al_n)A_n(z)-\textsf{v}_0'(z)\eqno(S_1)
$$
$$
1+(z-\al_n)(B_{n+1}(z)-B_n(z))=\bt_{n+1}A_{n+1}(z)-\bt_nA_{n-1}(z).
\eqno(S_2)
$$
We note here that $(S_1)$ and $(S_2)$ have been applied to random
matrix theory in \cite{TWFred}. It turns out that there is an
equation which gives better insight into the $\al_n$ and $\bt_n$, if
$(S_1)$ and $(S_2)$ are suitably combined. This equation was derived
in \cite{chen+its} for smooth weights and also holds at the presence
of discontinuities. Because its derivation is very simple, we recall
here the derivation.

Multiplying $(S_2)$ by $A_n(z)$, we see that the right hand side of
the resulting equation is a first order difference, while the left
hand side, with $(z-\al_n)A_n$ replaced by
$\quad\quad\quad\quad B_{n+1}(z)+B_{n}(z)+\textsf{v}_0'(z),$ is a first order difference
plus $A_n(z)$. Taking a telescope sum together with the appropriate
``initial condition''£¬
$$B_0(z)=A_{-1}(z)=0,$$ produces,
$$
B^2_n(z)+\textsf{v}_0'(z)B_n(z)+\sum_{j=0}^{n-1}A_j(z)=\bt_nA_n(z)A_{n-1}(z).
\eqno(S'_2).
$$
This last equation will be highly useful in what follows. The
equations $(S_1),$ $(S_2)$ and $(S_2')$ were also stated in
\cite{Mag}.

If $w_0$ is modified by the multiplication of "singular" factors
such as $|x-t|^{a}$, then the ladder operator relations, $(S_1),$
$(S_2)$ and $(S_2')$ remain valid; see \cite{chen+f}.

Let $\Psi(z)=P_n(z)$, eliminating $P_{n-1}(z)$ from the raising and
lowering operators gives \bea \label{2order.ode}
\Psi''(z)-\left(\textsf{v}^{\pr}(z)+\frac{A_n'(z)}{A_n(z)}\right)\Psi'(z)
+\left(B_n'(z)-B_n(z)\frac{A_n'(z)}{A_n(z)}+\sum_{j=0}^{n-1}A_j(z)\right)\Psi(z)=0,
\eea where we have used $(S_2')$ to simplify the coefficient of
$\Psi(z)$ in (\ref{2order.ode}).

For the problem at hand, $w_0(x)=x^{\al}(1-x)^{\bt},$ we must
suppose that $\al > 0$ and $\bt >0$ so that our weight is suitably
well-behaved. Therefore,
\begin{align}
w_0(x)&:=x^{\al}(1-x)^{\bt},\qquad x\in[0,1],\\
\textsf{v}_0(z)&:=-\al\log z-\bt\log(1-z),\non \\
\textsf{v}'_0(z)&=-\frac{\al}{z}-\frac{\bt}{z-1},\non \\
\frac{\textsf{v}_0'(z)-\textsf{v}_0'(y)}{z-y}&=\frac{\al}{zy}+\frac{\bt}{(z-1)(y-1)}.
\end{align}
Substituting (3.7) into the definitions of $A_n(z)$ and $B_n(z)$, followed by
integration by parts gives,
\begin{align}
A_n(z)&=\frac{R_n(t)}{z-t}-\frac{x_n(t)}{z-1}+\frac{x_n(t)-R_n(t)}{z},\non\\
B_n(z)&=\frac{r_n(t)}{z-t}-\frac{y_n(t)}{z-1}+\frac{y_n(t)-r_n(t)-n}{z},\non
\end{align}
where
\begin{align}
x_n(t)&:=\frac{\bt}{h_n}\int_{0}^{1}\frac{P_n^2(y)}{1-y}\;y^{\al}(1-y)^{\bt}w_{\textsf{J}}(y;t)dy,
\non\\
y_n(t)&:=\frac{\bt}{h_{n-1}}\int_{0}^{1}\frac{P_n(y)P_{n-1}(y)}{1-y}\;
y^{\al}(1-y)^{\bt}w_{\textsf{J}}(y;t)dy.\non
\end{align}

Note that, \begin{align}
A_n(z)&=O\left(\frac{1}{z^2}\right),\nn\\
B_n(z)&=-\frac{n}{z}+O\left(\frac{1}{z^2}\right),\nn
\end{align}
as $z\to\infty$, which is to be expected.
\vskip .35cm
We now begin with the characterization of the auxiliary quantities $R_n,\:r_n,\:x_n$
and $y_n.$

Substituting the expressions for $A_n(z)$ and $B_n(z)$ into $(S_1)$
and $(S_2')$, which are identities in $z$, and equating the residues
of the poles at $z=0, 1$ and $z=t$, we find 7 distinct difference
equations and 3 formulas which importantly performs the sums
$\sum_{j=0}^{n-1}R_j,$ $\sum_{j=0}^{n-1}x_j$ and
$\sum_{j=0}^{n-1}(R_j-x_j)$. \vskip .3cm From $(S_1),$ we find,
\begin{equation}
\label{ex1} r_{n+1}+r_n = (t-\al_n)R_n,
\end{equation}
\begin{equation} \label{ex2} -(y_{n+1}+y_n)=(\al_n-1)x_n+\bt,
\end{equation}
\begin{equation}
\label{ex3} y_{n+1}+y_n -r_{n+1}-r_n=2n+1+\al-\al_n(x_n-R_n).
\end{equation}

From $(S_2'),$ we find,
\begin{equation} \label{ex4}
\bt_nR_nR_{n-1}=r_n^2,
\end{equation}
\begin{equation} \label{ex5}
\bt_nx_nx_{n-1}=y_n^2+\bt y_n,
\end{equation}
\begin{equation}\label{ex6}
\bt_n(x_n-R_n)(x_{n-1}-R_{n-1})=(y_n-r_n-n)^2-\al(y_n-r_n-n),
 \end{equation}
\begin{equation}\label{ex7}
\bt_n(x_nR_{n-1}+x_{n-1}R_n)=(2n+\al+\bt)y_n-(2n+\al)r_n+2y_nr_n-n(n+\al),
\end{equation}

\begin{align}\label{ex8}
\sum_{j=0}^{n-1}R_j=&~\frac{(2n+\al+\bt)y_n-n(n+\al)}{t}\non\\
                   &+\frac{(2n+\al+\bt)(y_n-r_n)-n(n+\al)}{1-t},
\end{align}
\begin{align}\label{ex9}
\sum_{j=0}^{n-1}x_j=&~\frac{(2n+\al+\bt)(y_n-r_n)-n(n+\al)}{1-t}\non\\
                   &+(2n+\al+\bt)r_n+n(n+\al+\bt),
\end{align}
\begin{align}\label{ex10}
\sum_{j=0}^{n-1}(R_j-x_j)=&~\frac{(2n+\al+\bt)y_n-n(n+\al)}{t}\non\\
&-(2n+\al+\bt)r_n+n(n+\al+\bt).
\end{align}
Performing similar calculations on $(S_2)$, we find,
\begin{equation}\label{ex11}
(t-\al_n)(r_{n+1}-r_n)=\bt_{n+1}R_{n+1}-\bt_nR_{n-1},
\end{equation}
\begin{equation}\label{ex12}
(1-\al_n)(y_n-y_{n+1})=\bt_nx_{n-1}-\bt_{n+1}x_{n+1},
\end{equation}
\begin{equation}\label{ex13}
-\al_n(y_{n+1}-y_n+r_n-r_{n+1}-1)=\bt_{n+1}x_{n+1}
-\bt_nx_{n-1}+\bt_nR_{n-1}-\bt_{n+1}R_{n+1}.
\end{equation}

\noindent {\bf Remark 1.} It $t=0$, it is clear that
$r_n(0)=R_n(0)=0.$ So it follows from \eqref{ex7} that
$$y_n(0)=\frac{n(n+\al)}{2n+\al+\bt}.$$

\noindent {\bf Remark 2.} The equation \eqref{ex10} is a consistency
check on \eqref{ex8} and \eqref{ex9}.

\noindent {\bf Remark 3.} In \cite{basor+chen} and \cite{bacheh},
only $(S_1)$ and $(S_2')$ are required to provide a complete
description of the recurrence coefficient, while in \cite{chen+dai}
and \cite{dai+lun}, it was found that $(S_2)$ also plays a crucial
role. We shall see later it is also the case here.

We now manipulate the equations \eqref{ex1}--\eqref{ex13} with the
aim of expressing the recurrence coefficients $\al_n$ and $\bt_n$ in
terms of $r_n,\;R_n,\;x_n,\;y_n$ and of course $n$ and $t$. We start
with an easy lemma.
\begin{lemma}
We have the equations,
\begin{align}
\label{R_n and x_n}
x_n&=2n+1+\al+\bt+tR_n,\\
\al_n&=y_{n+1}-y_n+t(r_n-r_{n+1}), \label{aln y_n y_n+1}\\
\textsf{p}_1(n,t)&=-y_n(t)+tr_n(t) \label{p_1n-y_n-r_n}.
\end{align}
\end{lemma}
\vskip .3cm \noindent {\bf Proof.} The equation \eqref{R_n and x_n}
is obtained by summing \eqref{ex1}--\eqref{ex3}. If we substitute
the right hand side of \eqref{ex11} and \eqref{ex12} into
\eqref{ex13}, we find \eqref{aln y_n y_n+1}. The equation
\eqref{p_1n-y_n-r_n} follows from a telescopic sum of \eqref{aln y_n
y_n+1}, keeping in mind the ``initial'' conditions,
$r_0(t)=y_0(t)=0$, and $\sum_{j=0}^{n-1}\al_j=-\textsf{p}_1(n,t)$;
see \eqref{sum al_j}. \hfill $\Box$

%
\vskip .3cm Differentiating \eqref{p_1n-y_n-r_n} with respect to $t$
and taking into account \eqref{derivative of p_1n}, we find in the
next lemma
\begin{lemma}
\begin{equation}\label{dy and dr}
y_n'=tr_n'.
\end{equation}
\end{lemma}

Now, we are ready to express the recurrence coefficients $\al_n$ and
$\bt_n$ in terms of auxiliary quantities. In particular, we
emphasize here the importance of $\textsf{p}_1(n,t)$ as stated in
the following theorem.
\begin{theorem}
The recurrence coefficients $\al_n$ and $\bt_n$ have the following
representations:
\begin{align} \label{aln-r_n-y_n}
(2n+2+\al+\bt)\al_n&=2tr_n-2y_n-\bt+(2n+1+\al+\bt)t+(1-t)x_n,\nn\\
&=2\textsf{p}_{1}(n,t)-\bt+(2n+1+\al+\bt)t+(1-t)x_n.
\end{align}
and
\begin{align}
\label{equ for beta_n}
&(2n-1+\alpha+\beta)(2n+1+\alpha+\beta)\beta_n \nn\\
&=\big(y_n-tr_n\big)^2 +(2n+\al)tr_n+\big[\bt-(2n+\alpha+\beta)t\big]y_n+n(n+\al)t\nn\\
&=[\textsf{p}_1(n,t)]^2+(2n+\al)\textsf{p}_1(n,t)+(2n+\al+\bt)(1-t)y_n+n(n+\al)t\nn\\
&=[\textsf{p}_1(n,t)]^2+[-\bt+(2n+\al+\bt)t]\textsf{p}_1(n,t)\nn\\
&~~~+(2n+\al+\bt)t(1-t)\textsf{p}_1'(n,t)+n(n+\al)t.
\end{align}
\end{theorem}

\noindent \textbf{Proof:} First we eliminate $r_{n+1}$ from
\eqref{derivative of aln } and \eqref{ex1} to find,
\begin{equation}\label{r_n-R_n-al_nprime}
2r_n=(t-\al_n)R_n+\al_n'.
\end{equation}
In view of \eqref{derivative of aln }, we can replace $r_n-r_{n+1}$
in \eqref{aln y_n y_n+1} by $\al'_n$, i.e.,
\begin{equation}
\al_n=y_{n+1}-y_n+t\al_n'.\nn
\end{equation}
Now eliminate $y_{n+1}$ from above and \eqref{ex2}, we obtain
\begin{equation}\label{y_n-al_n-x_n}
2y_n=-\al_n+t\al_n'-\bt+(1-\al_n)x_n.
\end{equation}
The first equality in \eqref{aln-r_n-y_n} then follows from further
eliminating $\al_n'(t)$ by \eqref{r_n-R_n-al_nprime} and
\eqref{y_n-al_n-x_n} and the second equality is immediate on account
of \eqref{p_1n-y_n-r_n}.

To derive \eqref{equ for beta_n}, we use \eqref{R_n and x_n} to
eliminate $x_{n}$ and $x_{n-1}$ in \eqref{ex5}, and make use of
\eqref{ex4} to find,
\begin{align}\label{bn-1}
&\beta_n[c_nc_{n-1}+tc_nR_{n-1}+tc_{n-1}R_n
+t^2R_{n-1}R_{n}] \nonumber \\
&=\beta_n[c_nc_{n-1}+tc_nR_{n-1}+tc_{n-1}R_n]
+t^2r_n^2 \nn \\
&=y_n^2+\beta y_n.
\end{align}
where $c_n:=2n+1+\al+\bt$.

On the other hand, it follows from \eqref{ex7}, \eqref{R_n and x_n}
and \eqref{ex4} that
\begin{align}
&\beta_n(x_nR_{n-1}+x_{n-1}R_n) \nn \\
&=\beta_n\big((c_n+tR_n)R_{n-1}+(c_{n-1}+tR_{n-1})R_n\big)
\nn \\
&=\beta_n(c_nR_{n-1}+c_{n-1}R_n)+2tr_n^2 \nn \\
&= (2n+\al+\bt)y_n-(2n+\al)r_n+2y_nr_n-n(n+\al),
\end{align}
or equivalently,
\begin{align}\label{bn-2}
&\beta_n(c_nR_{n-1}+c_{n-1}R_n) \nn \\
&=(2n+\al+\bt)y_n-(2n+\al)r_n+2y_nr_n-2tr_n^2-n(n+\al).
\end{align}
Substitution of the above formula into \eqref{bn-1} gives us the
first equality of \eqref{equ for beta_n}. The last two equalities
follow from \eqref{p_1n-y_n-r_n} and \eqref{derivative of p_1n}.
\hfill $\Box$

Recall that $H_n(t)$ is the logarithmic derivative of $D_n(t)$
defined in \eqref{Hn}. The next theorem contains an important
identity expressing $H_n(t)$ in terms of the auxiliary quantities
$y_n$ and $r_n.$ This will in turn give us expressions of $y_n$ and
$r_n$ in terms of $H_n(t)$ and $H_n'(t).$
\begin{theorem}\label{H_n expression}
The logarithmic derivative of the Hankel determinant \eqref{Hn} can
be expressed in terms of $y_n$ and $r_n$ as
\begin{align}\label{Hn-y_n-r_n}
H_n(t)&=(2n+\al+\bt)(y_n-tr_n)-n(n+\al)\nn\\
&=-(2n+\al+\bt)\textsf{p}_1(n,t)-n(n+\al),
\end{align}
while
\begin{align}
r_n(t)&=-\frac{H_n'(t)}{2n+\al+\bt} \label{r_n-H_n}\\
y_n(t)&=\frac{-tH_n'(t)+H_n(t)+n(n+\al)}{2n+\al+\bt}
\label{y_n-H_n}.
\end{align}
\end{theorem}
\vskip .3cm {\bf Proof:} Combining \eqref{d log D_n}, which gives
the logarithmic derivative of $D_n(t)$ as the sum $\sum_jR_j$, and
\eqref{ex8}, which gives the sum in terms of $y_n$ and $r_n,$
leaves, \bea \frac{d}{dt}\ln
D_n=-\frac{(2n+\al+\bt)(y_n-r_n)-n(n+\al)}{1-t}
-\frac{(2n+\al+\bt)y_n-n(n+\al)}{t}.\nn \eea After a little
simplification, equation \eqref{Hn-y_n-r_n} follows. Taking a
derivative with respect to $t$ on the second equality in
\eqref{Hn-y_n-r_n}, it follows from \eqref{derivative of p_1n} that
$$
H_n'(t)=-(2n+\al+\bt)r_n(t),
$$
which is \eqref{r_n-H_n}. The derivation of \eqref{y_n-H_n} is now
obvious. \hfill $\Box$

\vskip .3cm \noindent {\bf Remark 4.} With reference to
\eqref{p_1n-y_n-r_n} and {\bf Remark 1}, we see at $t=0,$
$$
\textsf{p}_1(n,0)=-\frac{n(n+\al)}{2n+\al+\bt}.
$$
Integrating (3.33), we find,
$$
D_n(t)=D_n(0)\:\exp\left[(2n+\al+\bt)\int_{0}^{t}\:
\frac{\textsf{p}_1(n,s)-\textsf{p}_1(n,0)}{s(1-s)}ds\right],
$$
where $D_n(0)$ is Hankel determinant associated with ``unpertured'' Jacobi
weight
$$
(A+B)x^{\al}(1-x)^{\bt},
$$
and it is explicitly known.


Following what was done in \cite{basor+chen}, \cite{chen+dai},
\cite{dai+lun} and \cite{bacheh}, we derive two alternative
representations of $R_n$ in terms of $y_n$ and $r_n$, which is
stated in the next proposition.
\begin{prop}
The auxiliary quantity $R_n$ has the following representations
\begin{align}
R_n(t)= &\ \frac{(2n+1+\al+\bt)[l(r_n,
y_n,t)-(1-t)y_n'(t)]} { 2k(r_n, y_n,t)}, \label{R-r-y3} \\
\frac{1}{R_n(t)}=& \ \frac{l(r_n, y_n,t)+(1-t)y_n'(t)} {2(2n+ 1
+\al+\bt)r_n^2}, \label{R-r-y4}
\end{align}
where
\begin{equation} \label{l-r-y}
l(r_n,
y_n,t):=(2n+\alpha+\beta)y_n-(2n+\al)r_n-2tr_n^2+2y_nr_n-n(n+\al)
\end{equation}
and
\begin{align}\label{k-r-y}
k(r_n, y_n, t):=& ~\big(y_n-tr_n\big)^2 +(2n+\al)tr_n
+\big(\bt-(2n+\alpha+\beta)t\big)y_n \nonumber \\
& \ \ \ +n(n+\al)t
\end{align}
\end{prop}

\noindent \textbf{Proof:} Using \eqref{ex4}, we eliminate $R_{n-1}$
in \eqref{bn-2} and obtain
\begin{equation*}
\begin{aligned}
&(2n+1+\al+\bt)\frac{r_n^2}{R_n}+(2n-1+\al+\bt)\bt_nR_n
\\
&=(2n+\al+\bt)y_n-(2n+\al)r_n+2y_nr_n-2tr_n^2-n(n+\al).
\end{aligned}
\end{equation*}
Replacing $\bt_n$ in the above formula with the aid of \eqref{equ
for beta_n}, we have
\begin{align}
&\frac{(2n+1+\al+\bt)r_n^2}{R_n} + \frac{k(r_n,
y_n,t)}{2n+1+\al+\bt}R_n
\nonumber\\
&= (2n+\al+\bt+2r_n)y_n-(2n+\al)r_n-2tr_n^2-n(n+\al).
\label{R_n-r-y1}
\end{align}
This is a linear equation in $R_n$ and $1/R_n.$ We shall derive
another linear equation in $R_n$ and $1/R_n.$

To proceed we first note the first Toda equation \eqref{differential
equ for btn}, reproduced here;
\begin{equation}
\bt_n'=(R_{n-1}-R_n)\bt_n\nn
\end{equation}
Applying $\frac{d}{dt}$ to the first equality in \eqref{equ for
beta_n}, it follows from \eqref{ex4}, \eqref{equ for beta_n},
\eqref{dy and dr} and the above formula that
\begin{align}\label{R_n-r-y2}
&(2n-1+\al+\bt)(2n+1+\al+\bt)\frac{r_n^2}{R_n}-k(r_n,y_n,t)R_n \nn \\
&= 2tr_n^2+(2n+\al-2y_n)r_n-(2n+\alpha+\beta)y_n \nn \\
& \ \ \ +(2n+\alpha+\bt)(1-t)y_n'+n(n+\al).
\end{align}
The equations \eqref{R-r-y3} and \eqref{R-r-y4} now follow from
solving for $R_n$ and $1/R_n$ from (\ref{R_n-r-y1}) and
(\ref{R_n-r-y2}). \hfill $\Box$

\setcounter{equation}{0}
\section{Proof of Theorem \ref{thm of sigmaform} and derivation of Painlev\'e
VI}
\label{derivation of PVI}

With the above preparations, we are now ready to prove Theorem
\ref{thm of sigmaform} and show that $x_n(t)$ is, up to a linear
fractional transformation, a particular Painlev\'e VI.

\noindent \textit{Proof of Theorem \ref{thm of sigmaform}:}
Multiplying (\ref{R-r-y3}) and (\ref{R-r-y4}) gives us
\begin{equation} \label{r-r*}
(1-t)^2 [y_n'(t)]^2 = l^2(r_n, y_n,t) -  4k(r_n, y_n,t)r_n^2,
\end{equation}
where $l(r_n, y_n,t)$ and $k(r_n, y_n,t)$ are given by (\ref{l-r-y})
and (\ref{k-r-y}), respectively. Substituting \eqref{r_n-H_n} and
\eqref{y_n-H_n} into \eqref{r-r*} yields a second order non-linear
differential equation for $H_n$. Using (\ref{htil-def}) to replace
$H_n$ by $\sigma(t)$, we finally get (\ref{sigmaform}).

The initial conditions follow from \eqref{r_n-H_n} and the fact that
$r_n(0)=0$. This completes the proof of our theorem. \hfill $\Box$

\vskip .3cm The rest of this section is devoted to the derivation of
a differential equation satisfied by $x_n$. As this is going to
quite complicated algebraically, we will separate the proof into
several parts.
\begin{prop}
The auxiliary quantity $x_n$ has the following representations
\begin{align}
x_n(t)= &\ \frac{(2n+1+\al+\bt)[\widetilde l(r_n,
y_n,t)-t(1-t)y_n'(t)]} {2k(r_n, y_n,t)}, \label{x-r-y1} \\
\frac{1}{x_n(t)}=& \ \frac{\widetilde l(r_n, y_n,t)+t(1-t)y_n'(t)}
{2(2n+ 1 +\al+\bt)(\beta+y_n)y_n}, \label{x-r-y2}
\end{align}
where
\begin{equation} \label{tildel-r-y}
\widetilde l(r_n,
y_n,t):=2y_n^2+\big(2\beta-(2n+\alpha+\beta+2r_n)t\big)y_n
+(2n+\al)tr_n+n(n+\al)t
\end{equation}
and $k(r_n, y_n, t)$ is given in \eqref{k-r-y}.
\end{prop}

\noindent \textbf{Proof:} The equality \eqref{x-r-y1} follows from a
combination of \eqref{R_n and x_n} and \eqref{R-r-y3}. To show
\eqref{x-r-y2}, we use \eqref{R_n and x_n} to eliminate $R_n$ and
$R_{n-1}$ in \eqref{ex4} and get
\begin{equation}
\beta_nx_nx_{n-1}-\beta_nc_nx_{n-1}-\beta_nc_{n-1}x_n+\beta_nc_nc_{n-1}
=t^2r_n^2,
\end{equation}
where $c_n:=2n+1+\alpha+\beta$. On account of \eqref{ex5},
\eqref{equ for beta_n} and \eqref{k-r-y}, the above formulas is
equivalent to
\begin{equation}\label{equ x_n and 1/x_n}
y_n^2+\beta y_n-\frac{c_n(y_n^2+\beta y_n)}{x_n}-\frac{k(r_n, y_n,
t)x_n}{c_n}+k(r_n, y_n, t)=t^2r_n^2.
\end{equation}
Solving \eqref{equ x_n and 1/x_n} with the aid of \eqref{x-r-y1}
yields \eqref{x-r-y2}. \hfill $\Box$

Next, we express $r_n$ in terms of $x_n$, $x_n'$ and $y_n$.

Differentiating the first formula in \eqref{aln-r_n-y_n} with
respect to $t$, it follows from \eqref{derivative of aln } and
\eqref{dy and dr} that
\begin{equation}\label{r_n-r_n+1}
(2n+2+\al+\beta)(r_n-r_{n+1})=2r_n-x_n+(1-t)x_n'+2n+1+\al+\beta.
\end{equation}
By \eqref{ex1} and \eqref{R_n and x_n}, it is also easily seen that
\begin{equation}
(t-\alpha_n)(x_n-2n-1-\al-\beta)=t(r_n+r_{n+1}).
\end{equation}
Eliminating $r_{n+1}$ from the above two formulas, together with the
first equality of \eqref{aln-r_n-y_n}, we obtain an expression of
$r_n$ in terms of $x_n,\;x_n'$ and $y_n,$ in the next Lemma,
\begin{lemma}
\begin{equation}\label{rn-xn-yn}
r_n=-\frac{1}{2}+\frac{2n+1+\alpha+\beta+(1-t)x_n'}{2x_n}
-\frac{(2y_n+\beta-(1-t)x_n+t)(2n+1+\alpha+\beta-x_n)}{2tx_n}.
\end{equation}
\end{lemma}
As the final step involves rather large formulas, we will only give
a qualitative discussion here.

We begin by substituting \eqref{rn-xn-yn} into (\ref{x-r-y1}) and
(\ref{x-r-y2}). After some computations, we find a pair of
\emph{linear} equations in $y_n$ and $y_n'$, with coefficients
depending on $x_n,$ $x_n',\:n,\:t$. Solving this linear system gives
us
\begin{align}
y_n(t) &= F(x_n, x_n', n,t),\nn\\
y_n'(t)&= G(x_n, x_n', n,t),\nn
\end{align}
where the form of $F$ and $G$ are explicitly known. We note here
$x_n''$ \emph{does not} appear either in $F$ or in $G$. Since the
expressions are somewhat long, we are not prepared to write them
down.

From the obvious equality
$$
\frac{d}{dt} F(x_n, x_n',n,t) = y_n' (t)= G(x_n, x_n',n,t),
$$
it can be shown, after some computations, that the above equation is
equivalent to
\begin{align}\label{xn-ode} &\{(2n +\al+\bt)(2n
+\al+\bt
+1)+[(2n+\al+\bt+1)t-2(2n +\al+\bt)- 1]x_n(t)\nn\\
&-(t-1)x_n^2(t)+t(t-1)x_n'(t)\}\Phi(x_n,x_n', x_n'',n,t)=0,
\end{align}
where $\Phi$ is a functions that is explicitly known. Clearly, the
above formula yields two differential equations. One is a Riccati
equation, whose solutions are of no relevance to our problem and
should be discarded.

If we make the substitution \bea\label{wn-xn}
x_n(t)=(2n+1+\al+\bt)\frac{1-W_n(t)}{1-t} \eea in $\Phi=0$, we find
after some computations that, $W_n(t)$ satisfies the following
Painlev\'e VI equation
\begin{align}
W_n'' & = \frac{1}{2} \left( \frac{1}{W_n} + \frac{1}{W_n - 1} +
\frac{1}{W_n -t }\right) (W_n')^2 - \left( \frac{1}{t} + \frac{1}{t
- 1} + \frac{1}{W_n -t }\right) W_n' \nonumber \\
& \quad + \frac{W_n (W_n -1) (W_n - t)}{t^2 (t-1)^2} \left( a +
\frac{b\:t}{W_n^2} + \frac{c\:(t-1)}{(W_n -1)^2} + \frac{d\:
t(t-1)}{(W_n -t )^2} \right) \label{x-ode}
\end{align}
with
\begin{equation*}
a = \frac{(2n + \al + \bt + 1)^2}{2}, \quad b = -\frac{\al^2}{2},
\quad c = \frac{\bt^2}{2}, \quad d = \frac{1}{2}
\end{equation*}
and initial conditions
\begin{equation*}\label{initial condition}
W_n(0)=0, \qquad W_n'(0)=1,
\end{equation*}
as advertised in the introduction. The above initial conditions
follow from \eqref{R_n and x_n} and the fact $R_n(0)=0$.

\section{Proof of Theorem \ref{asy of gap probability}}
\label{proof of thm 2}

%

We begin with the integral representation of $D_n(t)$ given in
\eqref{integral representation}, which in the situation where $A=0$
and $B=1$ reads
\begin{align}
D_n(t)=\frac{1}{n!}\int_{[t,1]^n}[\Delta_n(x)]^{2}
\prod_{k=1}^{n}x_k^{\al}(1-x_k)^{\bt}dx_k,\quad
\De_n(x)=\prod_{1\leq i<j\leq n}(x_j-x_i).
\end{align}
A simple change of variable $x_k=y_k+t$ in the above multiple
integral yields
\begin{align}\label{D_n change of variable}
D_n(t)&=\frac{1}{n!}\int_{[0,1-t]^n}[\De_n(y)]^2\prod_{k=1}^{n}
(y_k+t)^{\al}(1-t-y_k)^{\bt}dy_k
\non\\
&=(1-t)^{n(n+\bt)}\frac{1}{n!}\int_{[0,1]^{n}}
[\De_n(y)]^2\:\prod_{k=1}^{n}[(1-t)y_k+t]^{\al}(1-y_k)^{\bt}dy_k
\non\\
&=(1-t)^{n(n+\bt)}\:t^{n\al}\frac{1}{n!}\int_{[0,1]^n}
\:\prod_{k=1}^{n}(1-y_k)^{\bt} \left[1-zy_k\right]^{\al}dy_k,
\end{align}
where $z=1-1/t$. As $t\to 1^{-},$ we find,
\begin{align}\label{D_n t to 1}
D_n(t)&\sim
(1-t)^{n(n+\bt)}\:\:\frac{1}{n!}\int_{[0,1]^n}[\De_n(y)]^2
\:\prod_{k=1}^{n}(1-y_k)^{\bt}dy_k.
\end{align}
To evaluate the multiple integral on the right hand side of
\eqref{D_n t to 1}, we set
\begin{align}
D_n[w_{\al,\bt}(-1,1)]
:=\frac{1}{n!}\int_{[-1,1]^n}[\De_n(x)]^2\prod_{k=1}^{n}
(1+x_k)^{\al}(1-x_k)^{\bt}dx_k.
\end{align}
Then, an elementary translation and scaling shows that
\begin{align}
D_n[w_{\al,\bt}(-1,1)]
&=2^{n(n+\al+\bt)}\:\frac{1}{n!}\int_{[0,1]^n}
[\De_n(y)]^2\prod_{k=1}^{n}y_k^{\al}(1-y_k)^{\bt}dy_k \non \\
&=2^{n(n+\al+\bt)}\:D_n[w_{\al,\bt}(0,1)]
\end{align}
By \cite[(1.6)]{basor+chen05}, we have
\begin{align}
D_n[w_{\al,\bt}(-1,1)]=2^{-n(n+\al+\bt)}\:(2\pi)^n\:F(\al,\bt)\:K(\al,\bt,n),
\end{align}
where $F(\al,\bt)$ and $K(\al,\bt,n)$ are given in \eqref{F(al,bt)}
and \eqref{K(al,bt)}, respectively. Hence, it follows from the above
two formulas that
\begin{align}\label{D_n(0,1)}
D_n[w_{\al,\bt}(0,1)]&=\frac{1}{n!}\int_{[0,1]^n}
[\De_n(y)]^2\prod_{k=1}^{n}y_k^{\al}(1-y_k)^{\bt}dy_k \non
\\
&=2^{-2n(n+\al+\bt)}\:(2\pi)^n\:F(\al,\bt)K(\al,\bt,n).
\end{align}
Note that $D_n(0)=D_n[w_{\al,\bt}(0,1)]$, combining \eqref{D_n
change of variable} and \eqref{D_n(0,1)} gives
\begin{align}
\frac{D_n(t)}{D_n(0)}
&\sim\frac{D_n[w_{0,\bt}(0,1)]}{D_n[w_{\al,\bt}(0,1)]}(1-t)^{n(n+\bt)}\non\\
&=2^{2n\al}\:\frac{F(0,\bt)K(0,\bt,n)}{F(\al,\bt)K(\al,\bt,n)}(1-t)^{n(n+\bt)},
\end{align}
which is \eqref{asy near 1}.

This finishes the proof of Theorem \ref{asy of gap
probability}.
\hfill $\Box$

\section*{Acknowledgement}
Lun Zhang is supported by FWO Flanders Project G.0427.09.

\end{document}